\documentclass[a4paper,12pt]{article}
\usepackage{amssymb, amsmath}
\usepackage[curve]{xypic}

\providecommand{\abs}[1]{\lvert#1\rvert}

\providecommand{\Ker}{\textnormal{Ker}}
\providecommand{\IIm}{\textnormal{Im}}

\providecommand{\Ker}{\textnormal{Ker}}

\providecommand{\rk}{\textnormal{rk}}
\providecommand{\Int}{\textnormal{Int}}

\providecommand{\PD}{\textnormal{PD}}
\providecommand{\cpt}{\textnormal{cpt}}

\providecommand{\chr}{\textnormal{char}}

\begin{document}

\begin{titlepage}
\titlepage
\rightline{SISSA 25/2009/FM}
\vskip 2.5cm
\centerline{ \bf \huge Gysin map and Atiyah-Hirzebruch }
\vskip 0.7cm
\centerline{ \bf \huge spectral sequence }
\vskip 1.7truecm

\begin{center}
{\bf \Large Fabio Ferrari Ruffino}
\vskip 1.5cm
\em 
International School for Advanced Studies (SISSA/ISAS) \\ 
Via Beirut 2, I-34014, Trieste, Italy\\
and Istituto Nazionale di Fisica Nucleare (INFN), sezione di Trieste

\vskip 2.5cm

\large \bf Abstract
\end{center}

\normalsize We discuss the relations between the Atiyah-Hirzebruch spectral sequence and the Gysin map for a multiplicative cohomology theory, on spaces having the homotopy type of a finite CW-complex. In particular, let us fix such a multiplicative cohomology theory $h^{*}$ and let us consider a smooth manifold $X$ of dimension $n$ and a compact submanifold $Y$ of dimension $p$, satisfying suitable hypotheses about orientability. We prove that, starting the Atiyah-Hirzebruch spectral sequence with the Poincar\'e dual of $Y$ in $X$, which, in our setting, is a simplicial cohomology class with coefficients in $h^{n-p}\{*\}$, if such a class survives until the last step, it is represented in $E^{n-p,\,0}_{\infty}$ by the image via the Gysin map of the unit cohomology class of $Y$. We then prove the analogous statement for a generic cohomology class on $Y$.

\vskip2cm

\vskip1.5\baselineskip

\vfill
 \hrule width 5.cm
\vskip 2.mm
{\small 
\noindent }
\begin{flushleft}
ferrariruffino@gmail.com
\end{flushleft}
\end{titlepage}

\newtheorem{Theorem}{Theorem}[section]
\newtheorem{Lemma}[Theorem]{Lemma}
\newtheorem{Corollary}[Theorem]{Corollary}
\newtheorem{Rmk}[Theorem]{Remark}
\newtheorem{Def}{Definition}[section]
\newtheorem{ThmDef}[Theorem]{Theorem - Defintion}


\section{Introduction}

Given a multiplicative cohomology theory, under suitable hypotheses we can define the Gysin map, which is a natural pushfoward in cohomology. Moreover, for a finite CW-complex or any space homotopically equivalent to it, we can construct the Atiyah-Hirzebruch spectral sequence, which relates the cellular cohomology with the fixed cohomology theory. In particular, the groups of the starting step of the spectral sequence $E^{p,\,q}_{1}(X)$ are canonically isomorphic to the groups of cellular cochains $C^{p}(X, h^{q}\{*\})$ for $\{*\}$ a fixed space with one point. Since the first coboundary $d^{p,\,q}_{1}$ coincides with the cellular coboundary, the groups $E^{p,\,q}_{2}(X)$ are canonically isomorphic to the cellular cohomology groups $H^{p}(X, h^{q}\{*\})$. The sequence stabilizes to $E^{p,\,q}_{\infty}(X)$ and, denoting by $X^{p}$ the $p$-skeleton of $X$, there is a canonical isomorphism:
\begin{equation}\label{EInftyIntro}
	E^{p,\,q}_{\infty}(X) \simeq \frac{\, \Ker ( h^{p+q}(X) \longrightarrow h^{p+q}(X^{p-1}) ) \,} {\Ker ( h^{p+q}(X) \longrightarrow h^{p+q}(X^{p}) )}
\end{equation}
i.e.\ $E^{p,\,q}_{\infty}$ can be described as the group of $(p+q)$-classes on $X$ which are $0$ when pulled back to $X^{p-1}$, up to classes which are $0$ when pulled back to $X^{p}$. Let us now consider an $n$-dimensional smooth manifold $X$ and a compact $p$-dimensional submanifold $Y$. For $i: Y \rightarrow X$ the embedding, we can define the Gysin map:
	\[i_{!}: h^{*}(Y) \longrightarrow \tilde{h}^{*+n-p}(X)
\]
which in particular gives a map $i_{!}: h^{0}(Y) \longrightarrow \tilde{h}^{n-p}(X)$. We assume that we have an oriented triangulation of $X$ restricting to a triangulation of $Y$ (this is always possible for $X$ orientable \cite{Munkres}): we require that $Y$ is a cycle in $C_{p}(X, h^{0}\{*\})$, identifying each simplex $\sigma$ of the triangulation with $\sigma \otimes_{\mathbb{Z}} 1$, for $1 \in h^{0}\{*\}$. Then, for $1 \in h^{0}(Y)$ defined as the pull-back of the unit $1 \in h^{0}\{*\}$ via the unique map $P: Y \rightarrow \{*\}$, we prove that $i_{!}(1)$ represents an element of $\Ker ( h^{p+q}(X) \rightarrow h^{p+q}(X^{p-1}) )$ (the latter being the numerator of \eqref{EInftyIntro}) and, if the Poincar\'e dual $\PD_{X}[Y] \in H^{n-p}(X, h^{0}\{*\})$ survives until the last step, its class in $E^{n-p,\,0}_{\infty}$ is represented exactly by $i_{!}(1)$. Similarly, for $\eta \in h^{0}\{*\}$, if the Poincar\'e dual of $Y \otimes \eta \in C_{p}(X, h^{0}\{*\})$ survives until $E^{n-p,\,0}_{\infty}$, its class is represented by $i_{!}(P^{*}\eta)$. More generally, without assuming $q = 0$, if $Y \otimes \alpha$ is a cycle in $C_{p}(X, h^{q}\{*\})$ for $\alpha \in h^{q}\{*\}$, and if $\PD_{X}[Y \otimes \alpha] \in H^{n-p}(X, h^{q}\{*\})$ survives until $E^{n-p,\,q}_{\infty}$, then its class in \eqref{EInftyIntro} is represented by $i_{!}(P^{*}\alpha)$. All the classes on $Y$ considered in these examples are pull-back of classes in $h^{*}\{*\}$: we will see that all the other classes give no more information.

The study of the relations between Gysin map and Atiyah-Hirzebruch spectral sequence was treated in \cite{FS} for K-theory, arising from the physical problem of relating two different classifications of D-brane charges in string theory. In this article we generalize the statement to any multiplicative cohomology theory.

The paper is organized as follows: in chapter \ref{SpecSeq} we briefly recall the basic theory of spectral sequences in order to show explicitely the maps needed in the following; in chapter \ref{OrientabGysinMap} we recall orientability, Thom isomorphism and Gysin map for a multiplicative cohomology theory; in chapter \ref{GysinAHSS} we state and prove the theorems providing the link between the Gysin map and the Atiyah-Hirzebruch spectral sequence.

\section{Spectral sequences}\label{SpecSeq}

\subsection{Review of Cartan-Eilenberg version}

We deal with spectral sequences in the axiomatic version described in \cite{CE}, chap.\ XV, par.\ 7, with the additional hypothesis of working with \emph{finite} sequences of groups. We also take into account the presence of the grading in cohomology. In particular, we suppose the following assignements are given for $p, p', p'' \in \mathbb{Z} \cup \{-\infty, +\infty\}$:
\begin{itemize}
	\item for $-\infty \leq p \leq p' \leq \infty$, abelian groups $H^{n}(p,p')$ for $n \in \mathbb{Z}$, such that $H^{n}(p, p') = H^{n}(0, p')$ for $p \leq 0$ and there exists $l \in \mathbb{N}$ such that $H^{n}(p, p') = H^{n}(p, +\infty)$ for $p' > l$ ($l$ does not depend on $n$ in our setting);
	\item for $p \leq p' \leq p''$, $a, b \geq 0$, $p+a \leq p'+b$, two maps:\footnote{The map $\delta$ is called in the same way in \cite{CE}. Instead, we introduce the name $\psi$ since the analogous map in \cite{CE} has no name.}
\begin{equation}\label{PsiDelta}
\begin{split}
	&\psi^{n}: H^{n}(p+a,p'+b) \rightarrow H^{n}(p,p')\\
	&\delta^{n}: H^{n}(p,p') \rightarrow H^{n+1}(p',p'')
\end{split}
\end{equation}
\end{itemize}
satisfying axioms (SP.1)-(SP.5) of \cite{CE}, p.\ 334. When the indices are not clear from the context, we also use the notations $(\psi^{n})^{p+a,p'+b}_{p,p'}$ and $(\delta^{n})^{p,p',p''}$ for the maps \eqref{PsiDelta}. We can describe the groups and the coboundaries of the spectral sequence in the following way:
\begin{equation}\label{Epr}
\begin{array}{ll}
	E^{p,\,q}_{r} = \IIm\bigl( H^{p+q}(p, p+r) \overset{\psi^{p+q}} \longrightarrow H^{p+q}(p-r+1, p+1) \bigr) & \textnormal{(\cite{CE}, formula (8) p.\ 318)} \\ \\
	d^{p,\,q}_{r} \,=\, (\delta^{p+q})^{p-r+1,p+1,p+r+1} \,\big\vert_{\IIm((\psi^{p+q})^{p,p+r}_{p-r+1,p+1})}\,:\\ \\
	\phantom{XXXXXXXXXXXXXXXXXX} E^{p,q}_{r} \longrightarrow E^{p+r,q-r+1}_{r} & \textnormal{(\cite{CE}, line 3 p.\ 319)} \\ \\
	F^{p,\,q}H = \IIm \bigl( H^{p+q}(p, +\infty) \overset{\psi^{p+q}}\longrightarrow H^{p+q}(0, +\infty) \bigr) & \textnormal{(\cite{CE}, line -10 p.\ 319)} \; .
\end{array}
\end{equation}
Then:
\begin{itemize}
	\item the groups $F^{p,\,q}H$ are a filtration of $H^{p+q}(0, +\infty)$;
	\item $\bigoplus_{p,q} E^{p,\,q}_{r+1} \simeq H\bigl( \bigoplus_{p,q} E^{p,\,q}_{r}, \bigoplus_{p,q} d^{p,\,q}_{r} \bigr)$ canonically, i.e.\ $E^{p,\,q}_{r+1} \simeq \Ker \, d^{p,\,q}_{r} / \IIm \, d^{p-r, \,q+r-1}_{r}$;
	\item the sequence $\{E^{p,\,q}_{r}\}_{r \in \mathbb{N}}$ stabilizes to $F^{p,\,q}H / F^{p+1,\,q-1}H$.
\end{itemize}
In particular, considering the following commutative diagram\footnote{The maps $\psi_{1}, \psi_{2}, \delta_{1}, \delta_{2}$ of the diagram are maps of the family \eqref{PsiDelta}; here and in the following we use this notation in order not to write too many indices.} (\cite{CE}, end of p.\ 318):
\begin{equation}\label{BoundaryDiagram}
\xymatrix{
H^{p+q}(p, p+r) \ar[r]^{\psi_{1} \quad\;\;} \ar[d]_{\delta_{1}} & H^{p+q}(p-r+1, p+1) \ar[d]^{\delta_{2}}\\
H^{p+q+1}(p+r, p+2r) \ar[r]^{\psi_{2} \;\;} & H^{p+q+1}(p+1,p+r+1)
}
\end{equation}
the following identities hold:
\begin{itemize}
	\item $\IIm(\psi_{1}) = E^{p,\,q}_{r}$ and $\IIm(\psi_{2}) = E^{p+r,\,q-r+1}_{r}$;
	\item $d^{p,\,q}_{r} = \delta_{2} \,\big\vert_{\IIm(\psi_{1})}\,:\, E^{p,\,q}_{r} \rightarrow E^{p+r,\,q-r+1}_{r}$.
\end{itemize}
The limit of the sequence $\bigoplus_{p} F^{p,\,q}H / F^{p+1,\,q-1}H$ can also be defined as (\cite{CE}, eq.\ (3) p.\ 316):
\begin{equation}\label{Limit}
	E^{p,\,q}_{0}H := E^{p,\,q}_{\infty} = \IIm \bigl( H^{p+q}(p, +\infty) \overset{\psi^{p+q}}\longrightarrow H^{p+q}(0, p+1) \bigr)
\end{equation}
i.e.\ $E^{p,\,q}_{0}H \simeq F^{p,\,q}H / F^{p+1,\,q-1}H$ canonically.

\subsection{Description of the isomorphisms}

We now explicitely show the isomorphisms and the maps we will need in the following. We postpone to the next subsection the proofs which cannot be found in \cite{CE}. Considering \eqref{BoundaryDiagram}, from the two diagrams:
\begin{equation}\label{BoundaryDiagram2}
\xymatrix{
H^{p+q-1}(p-r, p) \ar[r]^{\psi_{0}^{p+q} \qquad\;\;} \ar[d]^{\delta_{-1}^{p+q}} & H^{p+q-1}(p-2r+1, p-r+1) \ar[d]_{\delta_{0}^{p+q}}\\
H^{p+q}(p, p+r) \ar[r]^{\psi_{1}^{p+q} \qquad} \ar[d]^{\delta_{1}^{p+q}} & H^{p+q}(p-r+1, p+1) \ar[d]_{\delta_{2}^{p+q}}\\
H^{p+q+1}(p+r, p+2r) \ar[r]^{\psi_{2}^{p+q+1}\;\;} & H^{p+q+1}(p+1,p+r+1)\\
H^{p+q}(p, p+r+1) \ar[r]^{\psi_{3}^{p+q}} & H^{p+q}(p-r, p+1)
}
\end{equation}
we have that:
\begin{itemize}
	\item $\IIm(\psi_{1}^{p+q}) = E^{p,\,q}_{r}$;
	\item $d^{p,\,q}_{r} = \delta_{2}^{p+q} \,\big\vert_{\IIm(\psi_{1}^{p+q})}\,:\, E^{p,\,q}_{r} \rightarrow E^{p+r,\,q-r+1}_{r}$ and $d^{p-r,\,q+r-1}_{r} = \delta_{0}^{p+q} \,\big\vert_{\IIm(\psi_{0}^{p+q})}:\, E^{p-r,\,q+r-1}_{r} \rightarrow E^{p,\,q}_{r}$;
	\item $\IIm(\psi_{3}^{p+q}) = E^{p,\,q}_{r+1}$.
\end{itemize}
To find the isomorphism $E^{p,\,q}_{r+1} \simeq \Ker \, d^{p,\,q}_{r} \,/\, \IIm \, d^{p-r,\,q+r-1}_{r}$ we thus consider the map $\psi_{4}^{p+q}: H^{p+q}(p-r+1, p+1) \rightarrow H^{p+q}(p-r, p+1)$ which induces a surjection:
\begin{equation}\label{VarphiPqr}
	\varphi^{p,\,q}_{r} := \psi_{4}^{p+q}\vert_{\Ker \, d^{p,\,q}_{r}}: \Ker \, d^{p,\,q}_{r} \longrightarrow E^{p,\,q}_{r+1} \; .
\end{equation}
whose kernel is exactly $\IIm \, d^{p-r,\,q+r-1}_{r}$ (see subsection \ref{Proofs} for the proof).

Let us consider $E^{p,\,q}_{1} = H^{p+q}(p, p+1)$. Some elements lie in $\Ker \, d^{p,\,q}_{1}$, and they are mapped to $E^{p,\,q}_{2} \subset H^{p+q}(p-1, p+1)$ by $\varphi^{p,\,q}_{1}$, which is the restriction of $\psi^{p+q}: H^{p+q}(p, p+1) \rightarrow H^{p+q}(p-1, p+1)$ to such a kernel. We iterate the procedure: some elements of $E^{p,\,q}_{2}$ lie in $\Ker \, d^{p,\,q}_{2}$ and are mapped to $E^{p,\,q}_{3} \subset H^{p+q}(p-2, p+1)$ by $\varphi^{p,\,q}_{2}$, which is the restriction of $\psi^{p+q}: H^{p+q}(p-1, p+1) \rightarrow H^{p+q}(p-2, p+1)$. Thus, in the original group $E^{p,\,q}_{1} = H^{p+q}(p, p+1)$ we can consider the elements that survives to both these steps and we can map them directly to $E^{p,\,q}_{3} \subset H^{p+q}(p-2, p+1)$ via the composition $\psi^{p+q}: H^{p+q}(p, p+1) \rightarrow H^{p+q}(p-2, p+1)$. This procedure stops after $l$ steps (where $l$ is the number defined above such that $H^{n}(p, p') = H^{n}(p, +\infty)$ for any $p' > l$ and any $n$). In particular, we obtain a subset $A^{p,\,q} \subset E^{p,\,q}_{1}$ of \emph{surviving elements}, and a map:
\begin{equation}\label{MapAInfty}
	\varphi^{p,\,q}: A^{p,\,q} \subset E^{p,\,q}_{1} \longrightarrow E^{p,\,q}_{\infty}
\end{equation}
assigning to each surviving element its class in the last step. The map is simply the restriction of $\psi^{p+q}: H^{p+q}(p, p+1) \rightarrow H^{p+q}(0, p+1)$. The subgroup of surviving elements can be described as follows (see subsection \ref{Proofs} for the proof):
\begin{equation}\label{DescriptionApq}
	A^{p,\,q} = \IIm \bigl( H^{p+q}(p, +\infty) \overset{\psi_{5}^{p+q}}\longrightarrow H^{p+q}(p, p+1) \bigr)
\end{equation}
so that we can construct the commutative diagram:
\begin{equation}\label{DiagrE1EInfty}
\xymatrix{
	H^{p+q}(p, +\infty) \ar[rr]^{\psi_{7}^{p+q}} \ar[dr]_{\psi_{5}^{p+q}} & & H^{p+q}(0, p+1) \\
	& H^{p+q}(p, p+1) \ar[ur]_{\psi_{6}^{p+q}} &
}
\end{equation}
with $A^{p,\,q} = \IIm \, \psi_{5}^{p+q}$ and $\varphi^{p,\,q} = \psi_{6}^{p+q}\vert_{\IIm \, \psi_{5}^{p+q}}$.

\subsection{Proofs}\label{Proofs}

We now show the proofs of the statements of the previous subsection which cannot be found in \cite{CE}, at least in this axiomatic setting. The uninterested reader can skip to the next section.

Let us start with the map $\varphi^{p,\,q}_{r}$ defined in \eqref{VarphiPqr}. To prove that it is surjective, we consider the following commutative diagram:
\begin{tiny}
	\[\xymatrix{
	H^{p+q}(p, p+r+1) \ar[r]|{\psi_{5}} \ar@/^1pc/[rr]|{\psi_{8}} \ar@/^2pc/[rrrr]|{\psi_{3}} & H^{p+q}(p, p+r) \ar[r]|{\psi_{6}} \ar[d]|{\delta_{1}} \ar@/_1pc/[rr]|{\psi_{1}} & H^{p+q}(p, p+1) \ar[r]|(.43){\psi_{7}} \ar@/^1pc/[rr]|{\psi_{9}} \ar[dr]|{\delta_{3}} & H^{p+q}(p-r+1, p+1) \ar[r]|{\psi_{4}} \ar[d]|{\delta_{2}} & H^{p+q}(p-r, p+1) \\
	& H^{p+q+1}(p+r, p+2r) \ar[rr]|{\psi_{2}} & & H^{p+q+1}(p+1,p+r+1)
}\]
\end{tiny}
We show that:
\begin{itemize}
	\item \emph{The image of $\varphi^{p,\,q}_{r}$ is actually contained in $E^{p,\,q}_{r+1}$.} In fact, let us fix $a \in \Ker \, d^{p,\,q}_{r}$. Then there exists $b \in H^{p+q}(p, p+r)$ such that $\psi_{1}(b) = a$ and $\delta_{2} \circ \psi_{1}(b) = 0$. The latter is equivalent to $\delta_{3} \circ \psi_{6}(b) = 0$, i.e.\ $\psi_{6}(b) \in \Ker\, \delta_{3}$. By axiom (SP.4) p.\ 334 of \cite{CE} we have that $\Ker\, \delta_{3} = \IIm \, \psi_{8}$, thus there exists $c \in H^{p+q}(p, p+r+1)$ such that $\psi_{6}(b) = \psi_{8}(c)$, hence $\psi_{4}(a) = \psi_{4} \circ \psi_{1}(b) = \psi_{9} \circ \psi_{6}(b) = \psi_{9} \circ \psi_{8}(c) = \psi_{3}(c)$. This shows that the image under $\psi_{4}$ of $\Ker \, d^{p,\,q}_{r}$ is contained in the image of $\psi_{3}$.
	\item \emph{The image of $\varphi^{p,\,q}_{r}$ is the whole $E^{p,\,q}_{r+1}$.} Let us fix $a \in \IIm \, \psi_{3}$. Then there exists $b \in H^{p+q}(p, p+r+1)$ such that $a = \psi_{3}(b)$. Let us consider $c = \psi_{7} \circ \psi_{8}(b)$. Then $\psi_{4}(c) = a$ by construction, and we now show that $c \in \Ker \, d^{p,\,q}_{r}$. The fact that $c \in \IIm(\psi_{1})$ is obvious by construction, and $\delta_{2}(c) = \delta_{2} \circ \psi_{7} \circ \psi_{8}(b) = \delta_{3} \circ \psi_{8}(b) = 0$.
\end{itemize}
We now show that $\Ker \, \varphi^{p,\,q}_{r}$ is exactly $\IIm \, d^{p-r,\,q+r-1}_{r}$. We consider the following commutative diagram:
\begin{tiny}
	\[\xymatrix{
	H^{p+q-1}(p-r, p) \ar[r]|(.4){\psi_{10}} \ar[d]|{\delta_{-1}} \ar@/^1pc/[rr]|{\psi_{0}} & H^{p+q-1}(p-r, p-r+1) \ar[r]|(.45){\psi_{11}} \ar[dr]|{\delta_{4}} \ar[drrr]|{\delta_{5}} & H^{p+q-1}(p-2r+1, p-r+1) \ar[d]|{\delta_{0}} \\
	H^{p+q}(p, p+r) \ar[rr]|{\psi_{1}} & & H^{p+q}(p-r+1, p+1) \ar[r]|{\psi_{4}} \ar@/_1pc/[rr]|{\psi_{12}} & H^{p+q}(p-r, p+1) & H^{p+q}(p-r+1, p)
}\]
\end{tiny}
Let us consider $a \in \Ker \, \varphi^{p,\,q}_{r}$, which is equivalent to $a \in \Ker \, \psi_{4} \cap \Ker \, d^{p,\,q}_{r}$. Then, in particular, $a \in \Ker \, \psi_{4} \cap \IIm \, \psi_{1}$. We have to show that $a \in \IIm \, \delta_{0} \vert_{\IIm(\psi_{0})}$, i.e.\ that $a \in \IIm(\delta_{0} \circ \psi_{0})$. By axiom (SP.4) p.\ 334 of \cite{CE} we have that $\Ker \, \psi_{4} = \IIm \, \delta_{4}$, thus $a \in \IIm \, \delta_{4} \cap \IIm \, \psi_{1}$, so that there exist $b \in H^{p+q-1}(p-r, p-r+1)$ and $c \in H^{p+q}(p, p+r)$ such that $a = \delta_{4}(b) = \psi_{1}(c)$. Then $\psi_{12} \circ \delta_{4}(b) = \psi_{12} \circ \psi_{1}(c) = 0$, since $\psi_{12} \circ \psi_{1}: H^{p+q}(p, p+r) \rightarrow H^{p+q}(p-r+1, p)$ factorizes as $H^{p+q}(p, p+r) \rightarrow H^{p+q}(p, p) \rightarrow H^{p+q}(p-r+1, p)$ and $H^{p+q}(p, p) = 0$. Hence $\delta_{5}(b) = 0$, but $\Ker \, \delta_{5} = \IIm \, \psi_{10}$, thus there exists $d \in H^{p+q-1}(p-r, p)$ such that $b = \psi_{10}(d)$. This implies that $a = \delta_{0} \circ \psi_{0}(d)$ as claimed. Viceversa, let us show that for any $d \in H^{p+q-1}(p-r, p)$ it holds hat $\delta_{0} \circ \psi_{0}(d) \in \Ker \, \varphi^{p,\,q}_{r}$: in fact, $\psi_{4}(\delta_{0} \circ \psi_{0}(d)) = \psi_{4} \circ \delta_{4} \circ \psi_{10}(d) = 0$ since $\psi_{4} \circ \delta_{4} = 0$ by exactness.

\paragraph{}It remains to prove \eqref{DescriptionApq}. We show that the elements of $E^{p,\,q}_{1} = H^{p+q}(p, p+1)$ surviving until $E^{p,\,q}_{r}$ are:
	\[A^{p,\,q}_{r} = \IIm \bigl( H^{p+q}(p, p+r) \overset{\psi_{13}^{p+q}}\longrightarrow H^{p+q}(p, p+1) \bigr)
\]
from which \eqref{DescriptionApq} follows putting $r = +\infty$. We prove it by induction on $r$. For $r = 1$ the thesis is trivial. Let us consider the following diagram:
	\[\xymatrix{
	H^{p+q}(p, p+r+1) \ar[r]^{\quad \psi_{13}^{p+q}} & H^{p+q}(p, p+1) \ar[r]^{\psi_{14}^{p+q} \qquad} \ar[dr]_{\delta_{6}^{p+q}} & H^{p+q}(p-r+1, p+1) \ar[d]^{\delta_{2}^{p+q}} \\
	& & H^{p+q+1}(p+1, p+r+1).
}\]
As we said above, an element $a \in E^{p,\,q}_{1} = H^{p+q}(p, p+1)$ which survives until $E^{p,\,q}_{r}$ is mapped to $E^{p,\,q}_{r} \subset H^{p+q}(p-r+1, p+1)$ by $\psi_{14}^{p+q}$, and, if it survives also until $E^{p,\,q}_{r+1}$, its image lies in the kernel of $\delta_{2}^{p+q}$: thus $\delta_{2}^{p+q} \circ \psi_{14}^{p+q}(a) = 0$, which is equivalent to $\delta_{6}^{p+q}(a) = 0$. By exactness there exists $b \in H^{p+q}(p, p+r+1)$ such that $a = \psi_{13}^{p+q}(b)$, thus $a \in \IIm \, \psi_{13}^{p+q}$.

\subsection{Atiyah-Hirzebruch spectral sequence}

The Atiyah-Hirzebruch spectral sequence \cite{AH} relates the cellular cohomology of a finite CW-complex (or any space homotopically equivalent to it) to a generic cohomology theory $h^{*}$. For a finite simplicial complex $X$ we consider the natural filtration:
	\[\emptyset = X^{-1} \subset X^{0} \subset \cdots \subset X^{m} = X
\]
where $X^{i}$ is the $i$-th skeleton of $X$. The groups and maps of the spectral sequence are defined as follows:
\begin{itemize}
	\item $H^{n}(p,p') = h^{n}(X^{p'-1}, X^{p-1})$;
	\item $\psi^{n}: H^{n}(p+a, p'+b) \rightarrow H^{n}(p,p')$ is induced in cohomology by the map of couples $i: (X^{p'-1}, X^{p-1}) \rightarrow (X^{p'+b-1}, X^{p+a-1})$;
	\item $\delta^{n}: H^{n}(p,p') \rightarrow H^{n+1}(p',p'')$ is the composition of the map $\pi^{*}: h^{n}(X^{p'-1}, X^{p-1}) \rightarrow h^{n}(X^{p'-1})$ induced by the map of couples $\pi: (X^{p'-1}, \emptyset) \rightarrow (X^{p'-1}, X^{p-1})$, and the Bockstein map $\beta^{n}: h^{n}(X^{p'-1}) \longrightarrow h^{n+1}(X^{p''-1},X^{p'-1})$.
\end{itemize}
With these definitions all the axioms are satisfied, so that we can consider the corresponding spectral sequence $E^{p,\,q}_{r}(X)$. We briefly recall the structure of the first two and the last steps of such a sequence \cite{FS}. We have from \eqref{Epr} that $E^{p,\,q}_{1}(X) = H^{p+q}(p, p+1) = h^{p+q}(X^{p}, X^{p-1})$, thus $E^{p,\,q}_{1}(X) \simeq C^{p}(X, h^{q}\{*\})$ where $C^{p}(X, h^{q}\{*\})$ is the group of simplicial cochains with coefficients in $h^{q}\{*\}$ \cite{AH}. Moreover $d^{p,\,q}_{1}$ coincides with the simplicial coboundary operator, thus $E^{p,\,q}_{2}(X) \simeq H^{p}(X, h^{q}\{*\})$. For a more accurate review of the structure of cocycles and coboundaries we refer to \cite{FS}.

We denote by $i^{p}: X^{p} \rightarrow X$ the natural immersion and by $\pi^{p} : X \rightarrow X/X^{p}$ the natural projection for any $p$. For the Atiyah-Hirzebruch spectral sequence equation \eqref{Limit} is equivalent to:
\begin{equation}\label{EpinftyA}
	E^{p, \,q}_{\infty} = \IIm \bigl( \tilde{h}^{p+q}(X/X^{p-1}) \overset{\psi^{p+q}}\longrightarrow \tilde{h}^{p+q}(X^{p}) \bigr)
\end{equation}
where $\psi^{p+q}$ is the pull-back via $f^{p}: X^{p} \rightarrow X/X^{p-1}$ defined as $f^{p} = \pi^{p-1} \circ i^{p}$. Hence the following diagram commutes:\footnote{In the diagram we cannot say that $(i^{p})^{*} \circ (\pi^{p-1})^{*} = 0$ by exactness, since by exactness $(i^{p})^{*} \circ (\pi^{p})^{*} = 0$  at the same level $p$, as follows from $X^{p} \rightarrow X \rightarrow X/X^{p}$.}
\begin{equation}\label{EpinftyB}
\xymatrix{
\tilde{h}^{p+q}(X/X^{p-1}) \ar[dr]_{(\pi^{p-1})^{*}} \ar[rr]^{(f^{p})^{*}} & & \tilde{h}^{p+q}(X^{p})\\
& \tilde{h}^{p+q}(X) \ar[ur]_{(i^{p})^{*}}.
}
\end{equation}
The sequece $h^{p+q}(X, X^{p-1}) \overset{(\pi^{p-1})^{*}}\longrightarrow h^{p+q}(X) \overset{(i^{p-1})^{*}}\longrightarrow h^{p+q}(X^{p-1})$ is exact, i.e.\ $\IIm \, (\pi^{p-1})^{*} = \Ker \, (i^{p-1})^{*}$. Since trivially $\Ker \, (i^{p})^{*} \subset \Ker \,(i^{p-1})^{*}$, we obtain that $\Ker \, (i^{p})^{*} \subset \IIm \, (\pi^{p-1})^{*}$. Moreover:
	\[\IIm \, (f^{p})^{*} = \IIm \, \bigl( (i^{p})^{*} \circ (\pi^{p-1})^{*} \bigr) = \IIm \, \Bigl( (i^{p})^{*} \; \big\vert_{\IIm \, (\pi^{p-1})^{*}} \Bigr) \simeq \frac{\, \IIm \, (\pi^{p-1})^{*} \,} {\Ker \,(i^{p})^{*}} = \frac{\, \Ker \, (i^{p-1})^{*} \,} {\Ker \,(i^{p})^{*}}
\]
hence, finally:
\begin{equation}\label{EpqInfty}
	E^{p,\,q}_{\infty} \simeq \frac{\, \Ker \bigl( h^{p+q}(X) \longrightarrow h^{p+q}(X^{p-1}) \bigr) \,} {\Ker \bigl( h^{p+q}(X) \longrightarrow h^{p+q}(X^{p}) \bigr)}
\end{equation}
i.e.\ $E^{p,\,q}_{\infty}$ is made, up to canonical isomorphism, by $(p+q)$-classes on $X$ which are $0$ on $X^{p-1}$, up to classes which are $0$ on $X^{p}$.

\subsection{From the first to the last step}

We now see how to link the first and the last step of the sequence. In the diagram \eqref{DiagrE1EInfty} we know that an element $\alpha \in E^{p,\,q}_{1}$ survives until the last step if and only if $\alpha \in \IIm \, \psi_{5}^{p+q}$ and its class in $E^{p,\,q}_{\infty}$ is $\varphi^{p,\,q}(\alpha) = \psi_{6}^{p+q}(\alpha)$. We thus put, for $\alpha \in A^{p,\,q} = \IIm \, \psi_{5}^{p+q} \subset E^{p,\,q}_{1}$:
	\[\{\alpha\}_{E^{p,\,q}_{\infty}} := \varphi^{p,\,q}(\alpha) \; .
\]
For the Atiyah-Hirzebruch spectral sequence diagram \eqref{DiagrE1EInfty} becomes:
\begin{equation}\label{DiagrE1EInftyAHSS}
\xymatrix{
	\tilde{h}^{p+q}(X/X^{p-1}) \ar[rr]^{(f^{p})^{*}} \ar[dr]_{(i^{p,p-1})^{*}} & & \tilde{h}^{p+q}(X^{p}) \\
	& \tilde{h}^{p+q}(X^{p}/X^{p-1}) \ar[ur]_{(\pi^{p,p-1})^{*}} &
}
\end{equation}
for $\pi^{p, p-1}: X^{p} \rightarrow X^{p}/X^{p-1}$ the natural projection, $i^{p,p-1}: X^{p}/X^{p-1} \rightarrow X/X^{p-1}$ the natural immersion and $f^{p} = i^{p,p-1} \circ \pi^{p, p-1}$. Then $A^{p,\,q} = \IIm (i^{p,p-1})^{*}$ and $\varphi^{p,\,q} = (\pi^{p,p-1})^{*}\vert_{\IIm (i^{p,p-1})^{*}}$. Thus, the classes in $E^{p, \,q}_{1} = \tilde{h}^{p+q}(X^{p}/X^{p-1})$ surviving until the last step are the ones which are restrictions of a class defined on all $X / X^{p-1}$, and, for such a class $\alpha$:
\begin{equation}\label{FromOneToInfty}
	\{\alpha\}_{E^{p, \,q}_{\infty}} = (\pi^{p,p-1})^{*}(\alpha) \; .
\end{equation}

\section{Orientability and Gysin map}\label{OrientabGysinMap}

We consider the notion of multiplicative cohomology theory following \cite{Dold}. We recall that, if $h^{*}$ is a multiplicative cohomology theory, the coefficient group $h^{0}(\{*\})$, for $\{*\}$ a space with one point, is a commutative ring with unit. In fact, by the canonical homeomorphism $\{*\} \rightarrow \{*\} \times \{*\}$ we have a product $h^{0}(\{*\}) \times h^{0}(\{*\}) \rightarrow h^{0}(\{*\})$ which is associative. Moreover, skew-commutativity in this case coincides with commutativity, and $1$ is a unit also for this product.

Given a path-wise connected space $X$, we consider any map $p: \{*\} \rightarrow X$: by the path-wise connectedness of $X$ two such maps are homotopic, thus the pull-back $p^{*}: h^{*}(X) \rightarrow h^{*}(\{*\})$ is well defined.
\begin{Def}\label{DefRank} For $X$ a path-connected space we call \emph{rank} of a cohomology class $\alpha \in h^{n}(X)$ the class $\rk(\alpha) := (p^{*})^{n}(\alpha) \in h^{n}(\{*\})$ for any map $p: \{*\} \rightarrow X$.
\end{Def}
Let us consider the unique map $P: X \rightarrow \{*\}$.
\begin{Def}\label{TrivialClass} We call a cohomology class $\alpha \in h^{n}(X)$ \emph{trivial} if there exists $\beta \in h^{n}\{*\}$ such that $\alpha = (P^{*})^{n}(\beta)$. We denote by $1$ the class $(P^{*})^{0}(1)$.
\end{Def}
\begin{Lemma} For $X$ a path-wise connected space, a trivial chomology class $\alpha \in h^{n}(X)$ is the pull-back of its rank.
\end{Lemma}
\textbf{Proof:} Let $\alpha \in h^{n}(X)$ be trivial. Then $\alpha = (P^{*})^{n}(\beta)$ so that $\rk(\alpha) = (p^{*})^{n}(P^{*})^{n}(\beta) = (P \circ p)^{*\,n}(\beta) = \beta$, thus $\alpha = (P^{*})^{n}(\rk(\alpha))$. $\square$

\paragraph{}Let $\pi: E \rightarrow B$ be a fiber bundle with fiber $F$ and $E'$ a sub-bundle of $E$ with fiber $F' \subset F$. We have a natural diagonal map $\Delta_{\pi}: (E, E') \rightarrow (B \times E, B \times E')$ given by $\Delta_{\pi}(e) = (\pi(e), e)$ so that we can define the module structure:
\begin{equation}\label{Module2}
\begin{split}
	h^{i}(B) \times h^{j}(E, E') \overset{\times}\longrightarrow h^{i+j}(B \times E, B \times E') \overset{\Delta_{\pi}^{*}}\longrightarrow h^{i+j}(E, E') \; .
\end{split}
\end{equation}
\begin{Lemma}\label{Unitarity} The module structure \eqref{Module2} is unitary, i.e.\ $1 \cdot \alpha = \alpha$ for $1$ defined by \ref{TrivialClass}. More generally, for a trivial class $t = P^{*}(\eta)$, with $\eta \in h^{*}(\{*\})$, one has $t \cdot \alpha = \eta \cdot \alpha$.
\end{Lemma}
\textbf{Proof:} The thesis follows from the commutativity of the diagram:
\[\xymatrix{
	h^{i}(B) \times h^{j}(E, E') \ar[r]^{\times\quad} & h^{i+j}(B \times E, B \times E') \ar[r]^{\qquad\Delta_{\pi}^{*}} & h^{i+j}(E, E') \\
	h^{i}\{*\} \times h^{j}(E, E') \ar[r]^{\times\quad} \ar[u]^{(P^{*})^{i} \times 1^{j}} & h^{i+j}(\{*\} \times E, \{*\} \times E') \ar[ur]^{\simeq} \ar[u]^{((P \times 1)^{*})^{^{i+j}}} &
}\]
where the commutativity of the square follows directly from the naturality of the product, while the commutativity of triangle follows from the fact that $(P \times 1) \circ \Delta_{\pi}$ is the natural map $(E, E') \rightarrow (\{*\} \times E, \{*\} \times E')$ inducing the isomorphism $\simeq$. $\square$

\paragraph{}We now recall the notion of \emph{orientable} vector bundles with respect to a fixed \emph{multiplicative} cohomology theory. By hypothesis, there exists a unit $1 \in h^{0}(\{*\}) = \tilde{h}^{0}(S^{0})$. Since $S^{n}$ is homeomorphic to the $n$-th suspension of $S^{0}$, such a homeomorphism defines (by the suspension isomorphism) an element $\gamma^{n} \in \tilde{h}^{n}(S^{n})$ such that $\gamma^{n} = S^{n}(1)$ (clearly $\gamma^{n}$ is not the unit class since the latter does not belong to $\tilde{h}^{n}(S^{n})$). Moreover, given a vector bundle $E \rightarrow B$ with fiber $\mathbb{R}^{k}$, we have the canonical isomorphism in each fiber $F_{x} = \pi^{-1}(x)$:
\begin{equation}\label{IsoFF0FInfty}
h^{k}(F_{x}, (F_{x})_{0}) \simeq h^{k}(D^{k}_{x}, \partial D^{k}_{x}) \simeq h^{k}(D^{k}_{x} / \partial D^{k}_{x}, \partial D^{k}_{x} / \partial D^{k}_{x}) \simeq h^{k}(S^{k}, N)
\end{equation}
where the last isomorphism is non-canonical since it depends on the local chart ($N$ is the north pole of the sphere). However, since the homotopy type of a map from $S^{k}$ to $S^{k}$ is uniquely determined by its degree \cite{Hatcher} and a homeomorphism must have degree $\pm 1$, it follows that the last isomorphism of \eqref{IsoFF0FInfty} is canonical up to an overall sign, i.e.\ up to a multiplication by $-1$ in $h^{k}(S^{k}, N)$.
\begin{Def} Let $\pi: E \rightarrow B$ be a vector bundle of rank $k$ and $h^{*}$ a multiplicative cohomology theory in an admissible category $\mathcal{A}$ containing $\pi$. The bundle $E$ is called \emph{$h$-orientable} if there exists a class $u \in h^{k}(E, E_{0})$ such that for each fiber $F_{x} = \pi^{-1}(x)$ it satisfies $u\vert_{F_{x}} \simeq \pm\gamma^{k}$ under the isomorphism \eqref{IsoFF0FInfty}. The class $u$ is called \emph{orientation}.
\end{Def}

We now discuss some properties of $h$-orientations \cite{Rudyak}. The following lemma is very intuitive and can be probably deduced by a continuity argument; however, since we have not discussed topological properties of the cohomology groups, we give a proof not involving such problems. For a rank-$k$ vector bundle $\pi: E \rightarrow B$, let $(U_{\alpha}, \varphi_{\alpha})$ be a contractible local chart for $E$, with $\varphi_{\alpha}: \pi^{-1}(U_{\alpha}) \rightarrow U_{\alpha} \times \mathbb{R}^{k}$. Let us consider the compactification $\varphi_{\alpha}^{+}: \pi^{-1}(U_{\alpha})^{+} \rightarrow (U_{\alpha} \times \mathbb{R}^{k})^{+}$, restricting, for $x \in U_{\alpha}$, to $(\varphi_{\alpha})_{x}^{+}: E_{x}^{+} \rightarrow S^{k}$. Then we can consider the map:
\begin{equation}\label{VarphiHat}
	\hat{\varphi}_{\alpha,x} := ((\varphi_{\alpha})_{x}^{+\,-1})^{*\,k}: \tilde{h}^{k}(E_{x}^{+}) \longrightarrow \tilde{h}^{k}(S^{k}) \; .
\end{equation}
\begin{Lemma}\label{Continuity} Let $u$ be an $h$-orietation of a rank-$n$ vector bundle $\pi: E \rightarrow B$, let $(U_{\alpha}, \varphi_{\alpha})$ be a contractible local chart for $E$ and let $\hat{\varphi}_{\alpha,x}$ be defined by \eqref{VarphiHat}. Then $\hat{\varphi}_{\alpha,x}(u\vert_{E_{x}^{+}})$ is constant in $x$ with value $\gamma^{k}$ or $-\gamma^{k}$.
\end{Lemma}
\textbf{Proof:} Let us consider the map $(\varphi_{\alpha}^{+\,-1})^{*\,k}: \tilde{h}^{k}(\pi^{-1}(U_{\alpha})^{+}) \longrightarrow \tilde{h}^{k}((U_{\alpha} \times \mathbb{R}^{k})^{+})$ and let call $\xi := (\varphi_{\alpha}^{+\,-1})^{*\,k}(u\vert_{\pi^{-1}(U_{\alpha})^{+}})$. Since $(U_{\alpha} \times \mathbb{R}^{k})^{+} \simeq U_{\alpha} \times S^{k} \,/\, U_{\alpha} \times \{N\}$ canonically, we can consider the projection $\pi_{\alpha}: U_{\alpha} \times S^{k} \rightarrow U_{\alpha} \times S^{k} \,/\, U_{\alpha} \times \{N\}$. Then $\hat{\varphi}_{\alpha,x}(u\vert_{E_{x}^{+}}) = \xi\vert_{(\{x\} \times \mathbb{R}^{k})^{+}} \simeq \pi_{\alpha}^{*}(\xi)\vert_{\{x\} \times S^{k}}$. But, since $U_{\alpha}$ is contractible, the projection $\pi: U_{\alpha} \times S^{k} \rightarrow S^{k}$ induces an isomorphism in cohomology, so that $\pi_{\alpha}^{*}(\xi) = \pi^{*}(\eta)$ for $\eta \in h^{k}(S^{k})$, so that $\pi_{\alpha}^{*}(\xi)\vert_{\{x\} \times S^{k}} = \pi^{*}(\eta)\vert_{\{x\} \times S^{k}} \simeq \eta$, i.e.\ it is constant in $x$. By definition of orientation, its value must be $\pm\gamma^{k}$. $\square$

\begin{Theorem}\label{OrientedCharts} If a vector bundle $\pi: E \rightarrow B$ of rank $k$ is $h$-orientable, then given trivializing contractible charts $\{U_{\alpha}\}_{\alpha \in I}$ it is always possible to choose trivializations $\varphi_{\alpha}: \pi^{-1}(U_{\alpha}) \rightarrow U_{\alpha} \times \mathbb{R}^{k}$ such that $(\varphi_{\alpha}^{+})_{x}^{*\,k}(\gamma^{k}) = u\vert_{E_{x}^{+}}$. In particular, for $x \in U_{\alpha\beta}$ the homeomorphism $(\varphi_{\beta}\varphi_{\alpha}^{-1})^{+}_{x}: (\mathbb{R}^{k})^{+} \simeq S^{k} \longrightarrow (\mathbb{R}^{k})^{+} \simeq S^{k}$ satisfies $((\varphi_{\beta}\varphi_{\alpha}^{-1})^{+}_{x})^{*}(\gamma^{k}) = \gamma^{k}$.
\end{Theorem}
\textbf{Proof:} Choosen any local trivialization $\varphi_{\alpha}: \pi^{-1}(U_{\alpha}) \rightarrow U_{\alpha} \times \mathbb{R}^{k}$, it verifies $(\varphi_{\alpha}^{+})_{x}^{*\,k}(\gamma^{k}) = \pm u\vert_{E_{x}^{+}}$ by Lemma \ref{Continuity}. If the minus sign holds, it is enough compose $\varphi_{\alpha}$ to the pointwise reflection by an axes in $\mathbb{R}^{k}$, so that the compactified map has degree $-1$. $\square$

\begin{Def} An atlas satisfying the conditions of Theorem \ref{OrientedCharts} is called \emph{$h$-oriented atlas}.
\end{Def}

\begin{Lemma}\label{TwoOrientations} Let $\pi: E \rightarrow B$ be an $h^{*}$-orientable vector bundle of rank $k$, for $h^{*}$ a multiplicative cohomology theory. Then $E$ is orientable also with respect to the singular cohomology with coefficients in $h^{0}\{*\}$. Therefore, if $\chr (h^{0}\{*\}) > 2$, it is orientable in the usual sense. In particular, an atlas is $h$-oriented with respect to $u$ or $-u$ if and only if it is oriented.
\end{Lemma}
\textbf{Proof:} We call $\{\varphi_{\alpha\beta}\}$ the transition functions, and $\{\varphi_{\alpha\beta}^{+}\}$ their extension to the compactified fibers. Since $\varphi_{\alpha\beta}^{+}$ is a homeomorphism, it has degree $1$ or $-1$, and the degree of a map is independent of the cohomology theory \cite{Bredon}. If $\chr (h^{0}\{*\}) > 2$, an atlas is $h$-oriented, with respect to $u$ or $-u$, if and only if the degree of each $\varphi_{\alpha\beta}^{+}$ is $1$ and not $-1$, since $\varphi_{\alpha\beta}^{+}(\gamma^{k}) = \gamma^{k}$ (Theorem \ref{OrientedCharts}). The degree of $\varphi_{\alpha\beta}^{+}$ is $1$ if and only if the determinant of $\varphi_{\alpha\beta}$ is positive, thus the thesis follows. If $\chr (h^{0}\{*\}) = 2$ the thesis is trivial. $\square$

\paragraph{}Let $X$ be a compact smooth $n$-manifold and $Y \subset X$ a compact embedded $p$-dimensional submanifold such that the normal bundle $N(Y) = (TX\,\vert_{Y}) / \,TY$ is $h$-orientable. Then, since $Y$ is compact, there exists a tubular neighborhood $U$ of $Y$ in $X$ \cite{Bredon}, i.e.\ there exists an homeomorphism $\varphi_{U}: U \rightarrow N(Y)$.

If $i: Y \rightarrow X$ is the embedding, from this data we can naturally define an homomorphism, called \emph{Gysin map}:
	\[i_{!}: h^{*}(Y) \longrightarrow \tilde{h}^{*+n-p}(X) \; .
\]
In fact, we first apply the Thom isomorphism (\cite{Dold} page 7) $T: h^{*}(Y) \longrightarrow h^{*+n-p}_{\cpt}(N(Y)) = \tilde{h}^{*+n-p}(N(Y)^{+})$; then we naturally extend $\varphi_{U}$ to $\varphi_{U}^{+}: U^{+} \rightarrow N(Y)^{+}$ and apply $(\varphi_{U}^{+})^{*}: h^{*}_{\cpt}(N(Y)) \rightarrow h^{*}_{\cpt}(U)$; finally, considering the natural map $\psi: X \rightarrow U^{+}$ given by:
	\[\psi(x) = \left\{\begin{array}{ll}
	x & \text{if } x \in U \\
	\infty & \text{if } x \in X \setminus U
	\end{array}\right.
\]
we apply $\psi^{*}: \tilde{h}^{*}(U^{+}) \longrightarrow \tilde{h}^{*}(X)$. Summarizing:
\begin{equation}\label{GysinMap}
	i_{!}\,(\alpha) = \psi^{*} \circ \bigl(\varphi_{U}^{+}\bigr)^{*} \circ T \, (\alpha) \; .
\end{equation}

\paragraph{Remark:} One could try to use the immersion $i: U^{+} \rightarrow X^{+}$ and the retraction $r: X^{+} \rightarrow U^{+}$ to have a splitting $h(X) = h(U) \oplus h(X, U) = h(Y) \oplus h(X,U)$. But this is false, since the immersion $i: U^{+} \rightarrow X^{+}$ is not continuous: \emph{since $X$ is compact}, $\{\infty\} \subset X^{+}$ is open, but $i^{-1}(\{\infty\}) = \{\infty\}$, and $\{\infty\}$ is not open in $U^{+}$ since $U$ is non-compact.

\section{Gysin map and Atiyah-Hirzebruch spectral sequence}\label{GysinAHSS}

In this section we follow the same line of \cite{FS}, generalizing the discussion to any cohomology theory. We call $X$ a compact smooth $n$-dimensional manifold and $Y$ a compact embedded $p$-dimensional submanifold. We choose a finite triangulation of $X$ which restricts to a triangulation of $Y$ \cite{Munkres}. We use the following notation:
\begin{itemize}
	\item we denote the triangulation of $X$ by $\Delta = \{\Delta^{m}_{i}\}$, where $m$ is the dimension of the simplex and $i$ enumerates the $m$-simplices;
	\item we denote by $X_{\Delta}^{p}$ the $p$-skeleton of $X$ with respect to $\Delta$.
\end{itemize}
The same notation is used for other triangulations or simplicial decompositions of $X$ and $Y$. In the following theorem we need the definition of ``dual cell decomposition'' with respect to a triangulation: we refer to \cite{GH} pp.\ 53-54.
\begin{Theorem}\label{Triangulation} Let $X$ be an $n$-dimensional compact manifold and $Y \subset X$ a $p$-dimensional embedded compact submanifold. Let:
\begin{itemize}
	\item $\Delta = \{\Delta^{m}_{i}\}$ be a triangulation of $X$ which restricts to a triangulation $\Delta' = \{\Delta^{m}_{i'}\}$ of $Y$;
	\item $D = \{D^{n-m}_{i}\}$ be the dual decomposition of $X$ with respect to $\Delta$;
	\item $\tilde{D} \subset D$ be subset of $D$ made by the duals of the simplices in $\Delta'$.
\end{itemize}
Then, calling $\abs{\tilde{D}}$ the support of $\tilde{D}$:
\begin{itemize}
	\item the interior of $\abs{\tilde{D}}$ is a tubular neighborhood of $Y$ in $X$;
	\item the interior of $\abs{\tilde{D}}$ does not intersect $X_{D}^{n-p-1}$, i.e.:
	\[\abs{\tilde{D}} \cap X_{D}^{n-p-1} \subset \partial \abs{\tilde{D}} \; .
\]
\end{itemize}
\end{Theorem}
\textbf{Proof:} The $n$-simplices of $\tilde{D}$ are the duals of the vertices of $\Delta'$. Let $\tau = \{\tau^{m}_{j}\}$ be the first baricentric subdivision of $\Delta$ \cite{GH,Hatcher}. For each vertex $\Delta^{0}_{i'}$ in $Y$ (thought of as an element of $\Delta$), its dual is:
\begin{equation}\label{DTilde}
	\tilde{D}^{n}_{i'} = \bigcup_{\Delta^{0}_{i'} \in \tau^{n}_{j}} \tau^{n}_{j} \; .
\end{equation}
Moreover, if $\tau' = \{{\tau'}^{m}_{j'}\}$ is the first baricentric subdivision of $\Delta'$ (of course $\tau' \subset \tau$) and $D' = \{{D'}^{m}_{i'}\}$ is the dual of $\Delta'$ in $Y$, then (reminding that $p$ is the dimension of $Y$):
\begin{equation}\label{DPrimo}
	D'^{\,p}_{\;\,i'} = \bigcup_{\Delta^{0}_{i'} \in {\tau'}^{p}_{j'}} {\tau'}^{p}_{j'}
\end{equation}
and:
	\[\tilde{D}^{n}_{i'} \cap Y = D'^{\,p}_{\;\,i'} \; .
\]
Moreover, let us consider the $(n-p)$-simplices in $\tilde{D}$ contained in $\partial \tilde{D}^{n}_{i'}$ (for a fixed $i'$ in formula \eqref{DTilde}), i.e.\ $X^{n-p}_{\tilde{D}} \cap \tilde{D}^{n}_{i'}$: they intersect $Y$ transversally in the baricenters of each $p$-simplex of $\Delta'$ containing $\Delta^{0}_{i'}$: we call such baricenters $\{b_{1}, \ldots, b_{k}\}$ and the intersecting $(n-p)$-cells $\{\tilde{D}^{n-p}_{l}\}_{l = 1, \ldots, k}$. Since (for a fixed $i'$) $\tilde{D}^{n}_{i'}$ retracts on $\Delta^{0}_{i'}$, we can consider a local chart $(U_{i'}, \varphi_{i'})$, with $U_{i'} \subset \mathbb{R}^{n}$ neighborhood of $0$, such that:
\begin{itemize}
	\item $\varphi_{i'}^{-1}(U_{i'})$ is a neighborhood of $\tilde{D}^{n}_{i'}$;
	\item $\varphi_{i'}(D'^{\,p}_{\;\,i'}) \subset U_{i'} \cap (\{0\} \times \mathbb{R}^{p})$, for $0 \in \mathbb{R}^{n-p}$ (see eq. \eqref{DPrimo});
	\item $\varphi_{i'}(\tilde{D}^{n-p}_{l}) \subset U_{i'} \cap \bigl(\mathbb{R}^{n-p} \times \pi_{p}(\varphi_{i'}(b_{l}))\bigr)$, for $\pi_{p}: \mathbb{R}^{n} \rightarrow \{0\} \times \mathbb{R}^{p}$ the projection.
\end{itemize}
We now consider the natural foliation of $U_{i'}$ given by the intersection with the hyperplanes $\mathbb{R}^{n-p} \times \{x\}$ and its image via $\varphi_{i'}^{-1}$: in this way, we obtain a foliation of $\tilde{D}^{n}_{i'}$ transversal to $Y$. If we do this for any $i'$, by construction the various foliations glue on the intersections, since such intersections are given by the $(n-p)$-cells $\{\tilde{D}^{n-p}_{l}\}_{l = 1, \ldots, k}$, and the interior gives a $C^{0}$-tubular neighborhood of $Y$.

Moreover, a $(n-p-r)$-cell of $\tilde{D}$, for $r > 0$, cannot intersect $Y$ since it is contained in the boundary of a $(n-p)$-cell, and such cells intersect $Y$, which is done by $p$-cells, only in their interior points $b_{j}$. Being the simplicial decomposition finite, it follows that the interior of $\abs{\tilde{D}}$ does not intersect $X_{D}^{n-p-1}$. \\
$\square$

\paragraph{}We now consider quintuples $(X, Y, \Delta, D, \tilde{D})$ satisfying the following condition:
\begin{itemize}
	\item[$(\#)$] $X$ is an $n$-dimensional compact manifold and $Y \subset X$ a $p$-dimensional embedded compact submanifold such that $N(Y)$ is $h$-orientable. Moreover, $\Delta$, $D$ and $\tilde{D}$ are defined as in Theorem \ref{Triangulation}.
\end{itemize}

\begin{Lemma}\label{TrivialityXnp1} Let $(X,Y,\Delta,D,\tilde{D})$ be a quintuple satisfying $(\#)$, $U = \Int\abs{\tilde{D}}$ and $\alpha \in h^{*}(Y)$. Then:
\begin{itemize}
	\item there exists a neighborhood $V$ of $X \setminus U$ such that $i_{!}(\alpha)\vert_{V} = 0$;
	\item in particular, $i_{!}(\alpha) \,\vert_{X^{n-p-1}_{D}} = 0$.
\end{itemize}
\end{Lemma}
\textbf{Proof:} By equation \eqref{GysinMap}:
	\[i_{!}(\alpha) = \psi^{*} \beta \qquad \beta = (\varphi_{U}^{+})^{*} \circ \,T (\alpha) \in \tilde{h}^{*}(U^{+}) \; .
\]
Let $V_{\infty} \subset U^{+}$ be a contractible neighborhood of $\infty$, which exists since $U$ is a tubular neighborhood of a smooth manifold, and let $V = \psi^{-1}(V_{\infty})$. Then $\tilde{h}^{*}(V_{\infty}) \simeq \tilde{h}^{*}\{*\} = 0$, thus $\beta\vert_{V_{\infty}} = 0$ so that $(\psi^{*}\beta) \vert_{V} = 0$. By Theorem \ref{Triangulation} $X^{n-p-1}_{D}$ does not intersect the tubular neighborhood $\Int \abs{\tilde{D}}$ of $Y$, hence $X^{n-p-1}_{D} \subset V$, so that $(\psi^{*} \beta)\vert_{X^{n-p-1}_{D}} = 0$. $\square$

\subsection{Unit class}

We start by considering the case of the unit class $1 \in h^{0}(Y)$ (see def.\ \ref{TrivialClass}). Before we have assumed $X$ orientable for simplicity. We denote by $H$ the singular cohomology with coefficeints in $h^{0}\{*\}$: then the correct hypothesis is that $X$ must by $H$-orientable, since we need the Poincar\'e duality with respect to $H$. Therefore, the orientability of $X$ is necessary only if $\chr\,h^{0}\{*\} > 2$. If the normal bundle $N_{Y}X$ of $Y$ in $X$ is $h$-orientable, as in our hypotheses, then it is also $H$-orientable, thanks to Lemma \ref{TwoOrientations}. Actually, it also follows from the following argument. $Y$ is an $H$-orientable manifold: for $\chr\,h^{0}\{*\} = 2$ any bundle is orientable (thus also the tangent bundle $TY$), otherwise, being $Y$ a simplicial complex, in order to be a cycle in $C_{p}(X, h^{0}\{*\})$ it must be oriented as a simplicial complex, thus also as a manifold. Since also $X$ is $H$-orientable, it follows that both $TX\vert_{Y}$ and $TY$ are $H$-orientable, hence also $N_{Y}X$ is. Moreover, the atlas arising in the proof of Theorem \ref{Triangulation} is naturally $H$-oriented, as follows from the construction of the dual cell decomposition.

\begin{Theorem}\label{FirstTheorem} Let $(X,Y,\Delta,D,\tilde{D})$ be a quintuple satisfying $(\#)$ and $\Phi^{n-p}_{D}: C^{n-p}(X, h^{q}(\{*\}))$ $\rightarrow h^{n-p+q}(X_{D}^{n-p}, X_{D}^{n-p-1})$ be the standard canonical isomorphism. Let us define the natural projection and immersion:
	\[\pi^{n-p,\,n-p-1}: X_{D}^{n-p} \longrightarrow X_{D}^{n-p} / X_{D}^{n-p-1} \qquad\quad i^{n-p}: X_{D}^{n-p} \longrightarrow X
\]
and let $\PD_{\Delta}(Y)$ be the representative of $\PD_{X}[Y]$ given by the sum of the cells dual to the $p$-cells of $\Delta$ covering $Y$. Then:
	\[(i^{n-p})^{*}(i_{!}(1)) = (\pi^{n-p,\,n-p-1})^{*} ( \Phi_{D}^{n-p} (\PD_{\Delta}(Y) ) ) \; .
\]
\end{Theorem}
\textbf{Proof:} Let $U$ be the tubular neighborhood of $Y$ in $X$ stated in Theorem \ref{Triangulation}. We define the space $(U^{+})^{n-p}_{D}$ obtained considering the interior of the $(n-p)$-cells intersecting $Y$ transversally and compactifying this space to one point. The interiors of such cells forms exactly the intersection between the $(n-p)$-skeleton of $D$ and $U$, i.e.\ $X^{n-p}_{D} \, \vert_{U}$, since the only $(n-p)$-cells intersecting $U$ are the ones intersecting $Y$, and their interior is complitely contained in $U$, as stated in Theorem \ref{Triangulation}. If we close this space in $X$ we obtain the closed cells intersecting $Y$ transversally, whose boundary lies entirely in $X^{n-p-1}_{D}$. Thus the one-point compatification of the interior is:
	\[(U^{+})^{n-p}_{D} = \frac{\overline{X^{n-p}_{D} \, \vert_{U}}^{X}}{X^{n-p-1}_{D} \, \vert_{\partial U}}
\]
so that there is a natural inclusion $(U^{+})^{n-p}_{D} \subset U^{+}$ sending the denominator to $\infty$ (the numerator is exactly $X_{\tilde{D}}^{n-p}$ of Theorem \ref{Triangulation}). We also define:
	\[\psi^{n-p} = \psi \,\big\vert_{X^{n-p}_{D}}: \, X^{n-p}_{D} \longrightarrow (U^{+})^{n-p}_{D} \; .
\]
The latter is well-defined since the $(n-p)$-simplices outside $U$ and all the $(n-p-1)$-simplices are sent to $\infty$ by $\psi$. Calling $I$ the set of indices of the $(n-p)$-simplices in $D$, calling $S^{k}$ the $k$-dimensional sphere and denoting by $\dot{\cup}$ the one-point union of topological spaces, there are the following canonical homeomorphisms:
	\[\begin{split}
	& \xi^{n-p}_{X}: \pi^{n-p}(X_{D}^{n-p}) \overset{\simeq}\longrightarrow \dot{\bigcup_{i \in I}} \; S^{n-p}_{i} \\
	& \xi^{n-p}_{U^{+}}: \psi^{n-p}(X_{D}^{n-p}) \overset{\simeq}\longrightarrow \dot{\bigcup_{j \in J}} \; S^{n-p}_{j}
\end{split}\]
where $\{S^{n-p}_{j}\}_{j \in J}$, with $J \subset I$, is the set of $(n-p)$-spheres corresponding to the $(n-p)$-simplices with interior contanined in $U$, i.e.\ corresponding to $\pi^{n-p}\bigl(\overline{X^{n-p}_{D} \, \big\vert_{U}}\,\bigr)$. The homeomorphism $\xi^{n-p}_{U^{+}}$ is due to the fact that the boundary of the $(n-p)$-cells intersecting $U$ is contained in $\partial U$, hence it is sent to $\infty$ by $\psi^{n-p}$, while all the $(n-p)$-cells outside $U$ are sent to $\infty$: hence, the image of $\psi^{n-p}$ is homeomorphic to $\dot{\bigcup}_{j \in J} \; S^{n-p}_{j}$ sending $\infty$ to the attachment point. We define:
	\[\begin{split}
	\rho: \, \dot{\bigcup_{i \in I}} \; S^{n-p}_{i} \longrightarrow \dot{\bigcup_{j \in J}} \; S^{n-p}_{j}
\end{split}\]
as the natural projection, i.e.\ $\rho$ is the identity of $S^{n-p}_{j}$ for every $j \in J$ and sends all the spheres in $\{S^{n-p}_{i}\,\}_{i \in I \setminus J}$ to the attachment point. We have that:
	\[\xi^{n-p}_{U^{+}} \circ \psi^{n-p} = \rho \circ \xi^{n-p}_{X} \circ \pi^{n-p,\,n-p-1}
\]
hence:
\begin{equation}\label{PsiPiRho}
	(\psi^{n-p})^{*} \circ (\xi^{n-p}_{U^{+}})^{*} = (\pi^{n-p,\,n-p-1})^{*} \circ (\xi^{n-p}_{X})^{*} \circ \rho^{*} \; .
\end{equation}
We put $N = N(Y)$ and $\tilde{u}_{N} = (\varphi_{U}^{+})^{*}(u_{N})$, where $u_{N}$ is the Thom class of the normal bundle. By Lemma \ref{Unitarity} and equation \eqref{GysinMap} we have $i_{!}(1) = \psi^{*} \circ (\varphi_{U}^{+})^{*} (u_{N})$.
Then:
	\[(i^{n-p})^{*}(i_{!}(1)) = (i^{n-p})^{*} \psi^{*} (\tilde{u}_{N}) = (\psi^{n-p})^{*} \bigl( \tilde{u}_{N} \,\big\vert_{(U^{+})^{n-p}_{D}} \bigr)
\]
and
	\[(\xi^{n-p}_{X})^{*} \circ \rho^{*} \circ ((\xi^{n-p}_{U^{+}})^{-1})^{*} \bigl(\, \tilde{u}_{\mathcal{N}} \,\big\vert_{(U^{+})^{n-p}_{D}} \,\bigr) = \Phi_{D}^{n-p}(\PD_{\Delta}Y)
\]
since:
\begin{itemize}
	\item $\PD_{\Delta}(Y)$ is the sum of the $(n-p)$-cells intersecting $U$, oriented as the normal bundle;
	\item hence $((\xi^{n-p}_{X})^{-1})^{*} \circ \Phi_{D}^{n-p}(\PD_{\Delta}(Y))$ gives a $\gamma^{n-p}$ factor to each sphere $S^{n-p}_{j}$ for $j \in J$ and $0$ otherwise, orienting the sphere orthogonally to $Y$;
	\item but this is exactly $\rho^{*} \circ ((\xi^{n-p}_{U^{+}})^{-1})^{*} ( \tilde{u}_{N} \,\vert_{(U^{+})^{n-p}_{D}} )$ since by definition of orientability the restriction of $\tilde{\lambda}_{N}$ must be $\pm\gamma^{n}$ for each fiber of $N^{+}$. We must show that the sign ambiguity is fixed: this follows from the fact that the atlas arising from the tubular neighborhood in Theorem \ref{Triangulation} is $H$-oriented, as we pointed out at the beginning of this section. For the spheres outside $U$, that $\rho$ sends to $\infty$, we have that:
	\[\begin{split}
	\rho^{*} \bigl( \tilde{u}_{N} \,\big\vert_{(U^{+})^{n-p}_{D}} \bigr) \Big\vert_{\dot{\bigcup}_{i \in I \setminus J} \; S^{n-p}_{i}}
	&= \rho^{*} \bigl( \tilde{u}_{N} \,\big\vert_{\rho(\dot{\bigcup}_{i \in I \setminus J} \; S^{n-p}_{i})} \bigr)\\
	&= \rho^{*} \bigl( \tilde{u}_{N} \,\big\vert_{\{\infty\}} \bigr) = \rho^{*}(0) = 0 \; .
\end{split}\]
\end{itemize}
Hence, from equation \eqref{PsiPiRho}:
	\[\begin{split}
	i_{!}(Y \times \mathbb{C}) \, \big\vert_{X_{D}^{n-p}} &= (\psi^{n-p})^{*} \bigl( \tilde{u}_{N} \,\big\vert_{(U^{+})^{n-p}_{D}} \bigr)\\
	&= (\pi^{n-p,\,n-p-1})^{*} \circ (\xi^{n-p}_{X})^{*} \circ \rho^{*} \circ ((\xi^{n-p}_{U^{+}})^{-1})^{*}\bigl( \tilde{u}_{N} \,\big\vert_{(U^{+})^{n-p}_{D}} \bigr)\\
	&= (\pi^{n-p,\,n-p-1})^{*} \Phi_{D}^{n-p}(\PD_{\Delta}Y) \; .
\end{split}\]
$\square$

\paragraph{}Let us now consider any trivial class $P^{*}\eta \in h^{q}(Y)$. By Lemma \ref{Unitarity} we have that $P^{*}\eta \cdot u_{N} = \eta \cdot u_{N}$, hence Theorem \ref{FirstTheorem} becomes:
	\[(i^{n-p})^{*}(i_{!}(P^{*}\eta)) = (\pi^{n-p,\, n-p-1})^{*} ( \Phi_{D}^{n-p}( \PD_{\Delta}(Y \otimes \eta) )) \; .
\]
In fact, the same proof applies considering that $\eta \cdot u_{N}$ provides a factor $\eta \cdot \gamma^{n-p}$ instead of $\gamma^{n-p}$ for each spere of $N^{+}$, with $\eta \in h^{q}(\{*\}) \simeq \tilde{h}^{q}(S^{q})$.

\paragraph{}The following theorem encodes the link between Gysin map and AHSS.
\begin{Theorem}\label{SecondTheorem} Let $(X,Y,\Delta,D,\tilde{D})$ be a quintuple satisfying $(\#)$ and $\Phi^{n-p}_{D}: C^{n-p}(X, h^{q}(\{*\}))$ $\rightarrow h^{n-p+q}(X_{D}^{n-p}, X_{D}^{n-p-1})$ be the standard canonical isomorphism. Let us suppose that $\PD_{\Delta}Y$ is contained in the kernel of all the boundaries $d^{n-p, \,q}_{r}$ for $r \geq 1$. Then it defines a class:
	\[\{\Phi^{n-p}_{D}(\PD_{\Delta}(Y \otimes \eta))\}_{E^{n-p, \,q}_{\infty}} \in E^{n-p, \,q}_{\infty} \simeq \frac{\, \Ker ( h^{n-p+q}(X) \longrightarrow h^{n-p+q}(X^{n-p-1}) ) \,} {\Ker ( h^{n-p+q}(X) \longrightarrow h^{n-p+q}(X^{n-p}) )} \; .
\]
The following equality holds:
	\[\{\Phi^{n-p}_{D}(\PD_{\Delta}(Y \otimes \eta))\}_{E^{n-p, \,q}_{\infty}} = [i_{!}(P^{*}\eta)] \; .
\]
\end{Theorem}
\textbf{Proof:} By equations \eqref{EpinftyA} and \eqref{EpinftyB} we have:
\begin{equation}\label{Epinfty}
\xymatrix{
	E^{n-p, \,q}_{\infty} = \IIm \bigl( \tilde{h}^{n-p+q}(X/X_{D}^{n-p-1}) \ar[dr]_{(\pi^{n-p-1})^{*}} \ar[rr]^{\qquad\quad (f^{n-p})^{*}} & & \tilde{h}^{n-p+q}(X_{D}^{n-p}) \bigr)\\
& \tilde{h}^{n-p+q}(X) \ar[ur]_{(i^{n-p})^{*}}
}
\end{equation}
and, given a representative $\alpha \in \IIm \, (\pi_{n-r-1})^{*} = \Ker ( h^{n-p+q}(X) \longrightarrow h^{n-p+q}(X_{D}^{n-p-1}) )$, we have that $\{\alpha\}_{E^{n-p, \,q}_{\infty}} = (i^{n-p})^{*}(\alpha) = \alpha\,\vert_{X_{D}^{n-p}}$. Moreover, from \eqref{DiagrE1EInfty} we have the diagram:
\begin{equation}\label{DiagrE1EInfty2}
\xymatrix{
	E^{n-p, \,q}_{\infty} = \IIm \bigl( \tilde{h}^{n-p+q}(X/X^{n-p-1}_{D}) \ar[rr]^{\qquad\quad (f^{n-p})^{*}} \ar[dr]_{(i^{n-p,\,n-p-1})^{*}} & & \tilde{h}^{n-p+q}(X^{n-p}_{D}) \bigr)\\
	& \tilde{h}^{n-p+q}(X^{n-p}_{D} / X^{n-p-1}_{D}) \ar[ur]_{(\pi^{n-p,\,n-p-1})^{*}} \; . &
}
\end{equation}
where $i^{n-p,\,n-p-1}: X^{n-p}_{D} / X^{n-p-1}_{D} \rightarrow X/X^{n-p-1}$ is the natural immersion. We have that:
\begin{itemize}
	\item by formula \eqref{FromOneToInfty} the class $\{\Phi_{D}^{n-p} (\PD_{\Delta}(Y \otimes \eta))\}_{E^{n-p, \,q}_{\infty}}$ is given in diagram \eqref{DiagrE1EInfty2} by $(\pi^{n-p,\,n-p-1})^{*}(\Phi_{D}^{n-p}(\PD_{\Delta}(Y \otimes \eta)))$;
	\item by Lemma \ref{TrivialityXnp1} we have $i_{!}(1) \in \Ker ( h^{n-p+q}(X) \rightarrow h^{n-p+q}(X_{D}^{n-p-1}) )$, hence the class $[i_{!}(P^{*}\eta)]$ is well-defined in $E^{n-p, \,q}_{\infty}$, and, by exactness, $i_{!}(P^{*}\eta) \in \IIm \, (\pi^{n-p-1})^{*}$;
	\item by Theorem \ref{FirstTheorem} we have $(i^{n-p})^{*}( i_{!}(P^{*}\eta)) = (\pi^{n-p,\,n-p-1})^{*}(\Phi_{D}^{n-p}(\PD_{\Delta}(Y \otimes \eta)))$;
	\item hence $\{\Phi_{D}^{n-p} (\PD_{\Delta}(Y \otimes \eta))\}_{E^{n-p, \,q}_{\infty}} = [i_{!}(P^{*}\eta)]$.
\end{itemize}
$\square$

\begin{Corollary}\label{OrientableSurvives} Assuming the same data of the previous theorem, the fact that $Y$ has orientable normal bundle with respect to $h^{*}$ is a sufficient condition for $\PD_{\Delta}(Y)$ to survive until the last step of the spectral sequence. Thus, the Poincar\'e dual of any homology class $[\,Y\,] \in H_{p}(X, h^{q}\{*\})$ having a smooth representative with $h$-orientable normal bundle survives until the last step.
\end{Corollary}
\textbf{Proof:} we put together the diagrams \eqref{Epinfty} and \eqref{DiagrE1EInfty2}:
\begin{equation}\label{TwoDiagrams}
\xymatrix{
\tilde{h}^{n-p}(X/X^{n-p-1}_{D}) \ar[d]_{(i^{n-p,\, n-p-1})^{*}} \ar[rrr]^{\qquad (\pi^{n-p-1})^{*}} \ar[drrr]^{(f^{n-p})^{*}} & & & \tilde{h}^{n-p}(X) \ar[d]^{(i^{n-p})^{*}} \\
\tilde{h}^{n-p}(X^{n-p}_{D} / X^{n-p-1}_{D}) \ar[rrr]^{\qquad (\pi^{n-p,\,n-p-1})^{*}} & & & \tilde{h}^{n-p}(X^{n-p}_{D})
}
\end{equation}
and the diagram commutes being $\pi^{n-p,\,n-p-1} \circ i^{n-p,\, n-p-1} = i^{n-p} \circ \pi^{n-p-1}$. Under the hypotheses stated, we have that $i_{!}(1) \in \IIm (\pi^{n-p-1})^{*}$, so that $i_{!}(1) = (\pi^{n-p-1})^{*}(\alpha)$. Then $(i^{n-p})^{*}(\alpha) \in A^{n-p,\,0}$, so that it survives until the last step giving a class $(i^{n-p})^{*}(\pi^{n-p})^{*}(\alpha)$ in the last step. $\square$

\paragraph{}One could inquire if the condition of having $h$-orientable normal bundle is homology invariant. This is not true: let us consider the example of K-theory, for which a bundle is orientable if and only if it is a spin$^{c}$ bundle. In \cite{BHK} the authors show that in general, for a manifold $X$, there exist homologous submanifolds $Y$ and $Y'$, such that the normal bundle of $Y$ is spin$^{c}$, while the normal bundle of $Y'$ is not. Since the second step of the Atiyah-Hirzebruch spectral sequence coincides with the cohomology of $X$, this means that both $\PD_{\Delta}Y$ and $\PD_{\Delta'}Y'$ (for suitable $\Delta$ and $\Delta'$) survive until the last step, even if the normal bundle of $Y'$ is not orientable. Then, it is natural to inquire if it is true that a cohomology class survives until the last step if and only if it admits smooth representatives with orientable normal bundle, but we do not know the answer.

\subsection{Generic cohomology class}

If we consider a generic class $\alpha$ over $Y$ of rank $\rk(\alpha)$, we can prove that $i_{!}(E)$ and $i_{!}(P^{*}\rk(\alpha))$ have the same restriction to $X^{n-p}_{D}$: in fact, the Thom isomorphism gives $T(\alpha) = \alpha \cdot u_{N}$ and, if we restrict $\alpha \cdot u_{N}$ to a \emph{finite} family of fibers, which are transversal to $Y$, the contribution of $\alpha$ becomes trivial, so it has the same effect of the trivial class $P^{*}\rk(\alpha)$. We now prove this.

\begin{Lemma}\label{LineBundleXnp} Let $(X,Y,\Delta,D,\tilde{D})$ be a quintuple satisfying $(\#)$ and $\alpha \in h^{*}(Y)$ a class of rank $\rk(\alpha)$. Then:
	\[(i^{n-p})^{*}(i_{!}\alpha) = (i^{n-p})^{*}(i_{!}(P^{*}\rk\,\alpha)) \; .
\]
\end{Lemma}
\textbf{Proof:} Since $X^{n-p}_{D}$ intersects the tubular neighborhood in a finite number of cells corresponding under $\varphi_{U}^{+}$ to a finite number of fibers of the normal bundle $N$ attached to one point, it is sufficient to prove that, for any $y \in Y$, $(\alpha \cdot u_{N}) \, \vert_{N_{y}^{+}} = P^{*}\rk(\alpha) \cdot u_{N} \, \vert_{N_{y}^{+}}$. Let us consider the following diagram for $y \in B$:
\[\xymatrix{
	h^{i}(Y) \times h^{n}(N_{y}, N_{y}') \ar[r]^{\times \quad} \ar[d]_{(i^{*})^{i} \times (i^{*})^{n}} & h^{i+n}(Y \times N, Y \times N') \ar[d]^{(i\times i)^{*\,i+n}} \\
	h^{i}\{y\} \times h^{n}(N_{y}, N_{y}') \ar[r]^{\times \qquad} & h^{i+n}(\{y\} \times N_{y}, \{*\} \times N'_{y}) \; .
}\]
The diagram commutes by naturality of the product, thus $(\alpha \cdot u_{N}) \, \vert_{N_{y}^{+}} = \alpha\vert_{\{y\}} \cdot u_{N}\vert_{N_{y}^{+}}$. Thus, we just have to prove that $\alpha\vert_{\{y\}} = (P^{*}\rk(\alpha))\,\vert_{\{y\}}$, i.e. that $i^{*}\alpha = i^{*}P^{*}p^{*}\alpha = (p \circ P \circ i)^{*}\alpha$. This immediately follows from the fact that $p \circ P \circ i = i$.
$\square$

\paragraph{}In the previous theorems we started from the first step of the spectral sequence, therefore we had to choose a simplicial decomposition of $X$. Anyway, if we start from the second step, we loose the dependence on the triangulation \cite{AH}.


\begin{thebibliography}{99}

\bibitem{AH} M. Atiyah and F. Hirzebruch, \emph{Vector Bundles and Homogeneous Spaces}, Michael Atiyah: Collected works, v. 2
\bibitem{BHK} C. Bohr, B. Hanke, D. Kotschick, \emph{Cycles, submanifolds and structures on normal bundles}, Manuscripta Math. 108 (2002), 483--494
\bibitem{Bredon} G. E. Bredon, \emph{Topology and geometry}, Springer-Verlag, 1993
\bibitem{CE} H. Cartan and S. Eilenberg, \emph{Homological algebra}, Princeton University Press, 1956
\bibitem{Dold} A. Dold, \emph{Relations between ordinary and extraordinary homology}, Colloquium on Algebraic Topology, Institute of Mathematics Aarhus University, 1962, pp. 2-9
\bibitem{FS} F. Ferrari Ruffino and R. Savelli, \emph{Comparing two approaches to the K-theory classification of D-branes}, to appear in Journal of Geometry and Physics, arXiv:0805.1009
\bibitem{GH} P. Griffiths and J. Harris, \emph{Principles of algebraic geometry}, John Wiley \& Sons, 1978
\bibitem{Hatcher} A. Hatcher, \emph{Algebraic topology}, Cambridge university press, 2002
\bibitem{Munkres} J.R. Munkres, \emph{Elementary Differential Topology}, Princeton University Press, 1968
\bibitem{Rudyak} Y. B. Rudyak, \emph{On Thom spectra, orientability and cobordism}, Springer monographs in mathematics

\end{thebibliography}
\end{document}